\documentstyle[amssymb]{article} 
\def\eqalign#1{%
\null \,\vcenter {\openup \jot \ialign {\strut \hfil $\displaystyle {
##}$&$\displaystyle {{}##}$\hfil \crcr #1\crcr }}\,}

\newcommand{\dom}{{\rm dom}}

\newcommand{\QED}{\vrule width 6pt height 6pt depth 0pt \vspace{0.1in}}
\newcommand{\forces}{\mathrel{\|}\joinrel\mathrel{-}}
\newcommand{\F}{{\cal F}}

\newcommand{\V}{{\bf V}}
\def\<#1\>{\langle #1 \rangle}

\newtheorem{theorem}{Theorem}[section]
\newtheorem{lemma}[theorem]{Lemma}

\newcommand{\lesdot}{\mathrel{\mathord{<}\!\!\raise
0.8 pt\hbox{$\scriptstyle\circ$}}}
\newcommand{\Proof}{{\sc Proof} \hspace{0.2in}}

\newcommand{\lft}[2]{\mathopen\ifcase#1{}\oo\or
                        \big#2\or\Big#2\else\oo\fi}
\newcommand{\rgt}[2]{\mathclose\ifcase#1{}\oo\or
                        \big#2\or\Big#2\else\oo\fi}
\def\lh(#1){|#1|}

\begin{document}

\title{Intersection of $< 2^{\aleph_{0}}$ ultrafilters may have measure zero}
\author{Tomek Bartoszynski\thanks{The author thanks the  Lady Davis
Fellowship
Trust for full support} \\
Boise State University\\
Boise, Idaho\\
and\\
Hebrew University \\
Jerusalem
\and
Saharon Shelah\thanks{Partially supported by Basic Research Fund,
Israel Academy of Sciences, publication 436}\\
Hebrew University\\
Jerusalem}
\maketitle
\begin{abstract}
We show that it is consistent with ZFC that the intersection of some
family of less than
$2^{\aleph_{0}}$ ultrafilters can have measure zero.
This answers a question of D. Fremlin.
\end{abstract}
The goal of this paper is to prove the theorem in the title.
The solution is due to the second author.

Throughout the paper we use standard notation. All the filters are
assumed to be non-principal filters on $\omega$, i.e., they
contain the Frechet filter ${\cal F}_0 = \left\{X \subseteq \omega :
|\omega-X|<\aleph_{0}\right\}$.
We identify filters with subsets of
$2^\omega$. In this way the question about  measurability and the
Baire property of filters on $\omega$ makes sense.
Let $\mu$ denote the standard measure on $2^\omega$ and $\mu^\star$
and $\mu_\star$ the outer and inner measure respectively.

We would like to thank the referee for many helpful remarks which
improved the quality of this paper.
\section{Introduction}
\setcounter{theorem}{0}
In this section we present several results concerning intresections of
filters on $\omega$.
The following classical result is a starting point for all subsequent
theorems.

\begin{theorem}[Sierpinski]\label{sierpinski}
Suppose that ${\cal F}$ is a filter on $\omega$. Then ${\cal F}$ has
measure zero or is nonmeasurable. Similarly, it either is meager or
does not have the Baire property.
If ${\cal F}$ is an ultrafilter then ${\cal F}$ is nonmeasurable and
does not have the Baire property. $\QED$
\end{theorem}

In particular if ${\cal F}$ is a filter then $\mu_\star({\cal F})=0$
and if ${\cal F}$ is an ultrafilter then $\mu^\star({\cal F})=1$.

First consider the filters which do not have the Baire property.

\begin{theorem}[Talagrand {[T]}]\label{2.1}
\ \begin{enumerate}
\item The intersection of countably many filters without the 
Baire property is
a filter without the Baire property.
\item Assume ${\bf MA}$. Then the 
intersection of $< 2^{\aleph_{0}}$ filters
without the Baire property is
a filter without the Baire property. $\QED$
\end{enumerate}
\end{theorem}

For ultrafilters we have  a much stronger result.

\begin{theorem}[Plewik {[P]}]\label{plew}
The intersection of $< 2^{\aleph_{0}}$ ultrafilters is a filter which does
not have the
Baire property. $\QED$
\end{theorem}

For Lebesgue measure the situation is more complicated.

\begin{theorem}[Talagrand {[T]}]
Let $\{\F_n : n \in \omega\}$ be a family of nonmeasurable filters.
Then $\F=\bigcap_{n \in \omega} \F_n$ is nonmeasurable.
\end{theorem}

If we consider uncountable families of filters the analog of
\ref{2.1} is no longer true.

\begin{theorem}[Fremlin {[F]}]\label{fr}
Assume Martin's Axiom.

Then there exists a family
$\left\{{\cal F}_{\xi}~ :~ \xi~ <~ 2^{\aleph_{0}}\right\}$ of nonmeasurable filters
such that
$\bigcap_{\xi \in I} {\cal F}_{\xi}$ is a measurable filter
for every uncountable set $I \subset 2^{\aleph_{0}}$.
In particular there exists a family of $\aleph_{1}$ nonmeasurable
filters with measurable intersection.~$\QED$
\end{theorem}

The next theorem shows that the above pathology cannot happen if we assume
stronger measurability properties.

Let $\vec{p} = \<p_{n}:n \in \omega \>$ be a sequence of
reals such that $p_{n}~ \in~ (0,\frac{1}{2}]$ for all $n \in \omega$.

Define $\mu_{\vec{p}}$ to be the product measure on
$2^{\omega}$ such that for all $n$,
$$\mu_{\vec{p}}(\left\{x \in 2^{\omega}  :  x(n)=1\right\})  = 
p_{n}$$ and
$$\mu_{\vec{p}}(\left\{x  \in  2^{\omega}  :   x(n)=0\right\})  =  1 
-  p_{n}.$$
Notice that if $p_{n} = \frac{1}{2}$ for all $n$ then $\mu_{\vec{p}}$
is the usual measure on $2^{\omega}$.

\begin{theorem}[Bartoszynski {[Ba]}]\label{2.4}
Assume ${\bf MA}$. Let $\mu_{\vec{p}}$ be a
measure such that $\lim_{n \rightarrow \infty} p_{n}=0$ and let
$\left\{{\cal F}_{\xi} :
\xi < \lambda < 2^{\omega}\right\}$ be a family of $\mu_{\vec{p}}$-nonmeasurable
filters.
Then
$$\bigcap_{\xi < \lambda} {\cal F}_{\xi}
\hbox{ is a Lebesgue nonmeasurable filter} .\  \QED$$
\end{theorem}

Finally notice that the additional assumptions in \ref{2.1} and
\ref{fr} are necessary. Namely, we have the following result.
\begin{theorem}
It is consistent with ZFC that there exists a family of filters
${\cal A}$ of size $2^{\aleph_{0}}$ such that
\begin{enumerate}
\item ${\cal A}$ consists of filters which do not have the Baire property
and the intersection of any uncountable subfamily of ${\cal A}$ is
equal to ${\cal F}_0$,
\item ${\cal A}$ consists of non-measurable filters
and the intersection of any uncountable subfamily of ${\cal A}$ is
equal to ${\cal F}_0$.
\end{enumerate}
\end{theorem}
\Proof
1) Let $\V$ be a model of ZFC satisfying CH. Let $\<c^\xi_\eta : \xi <
\omega_1, \eta < \kappa\>$ be a generic sequence of Cohen reals over
$\V$.
Let ${\cal F}_\eta$ be the filter generated by $\left\{c^\xi_\eta : \xi <
\omega_1\right\}$ for $\eta < \kappa$.

One easily checks that this family of filters has the required properties.

2) Use a sequence of random  instead of Cohen reals. $\QED$.

\section{Intersection of ultrafilters}
\setcounter{theorem}{0}
In this section we show that the analog of \ref{plew} is not true.
\begin{theorem}\label{main}
It is consistent with ZFC that the intersection of some family of
$<2^{\aleph_{0}}$
ultrafilters has measure zero.
\end{theorem}
\Proof
We start with the following observation.
\begin{lemma}
Suppose that $\left\{{\cal F}_\xi : \xi < \kappa\right\}$ is a family of filters
on $\omega$ such that
$\mu_\star(\bigcup_{\xi < \kappa} {\cal F}_\xi)>0 \ .$
Then
$$\mu(\bigcap_{\xi < \kappa} {\cal F}_\xi)=0 \ .$$
Moreover, if $\F_\xi$'s are ultrafilters, $\bigcap_{\xi <
\kappa} \F_\xi$ has measure zero iff $\bigcup_{\xi < \kappa} \F_\xi$
has measure 1.
\end{lemma}
\Proof
Suppose that
$\mu_\star(\bigcup_{\xi < \kappa} {\cal F}_\xi)>0$. Since all filters
are assumed to be non-principal we know that
$\mu_\star(\bigcup_{\xi < \kappa} {\cal F}_\xi)=1$. Let $A \subseteq
2^\omega$ be a measure 1 set contained in
$\bigcup_{\xi < \kappa} {\cal F}_\xi$. Define $A^\star = \left\{\omega-X :
X \in A\right\}$. Clearly $\mu(A^\star)=1$.

We claim that
$$A^\star \cap \bigcap_{\xi < \kappa} \F_\xi = \emptyset \ .$$
Suppose that $X \in A^\star$. Then $\omega-X \in A$ and hence
$\omega-X \in \F_\eta$ for some $\eta < \kappa$.
It follows that $X \not \in \bigcap_{\xi < \kappa} \F_\xi$.

Conversely, suppose that $\bigcap_{\xi < \kappa} \F_\xi$ has measure
zero. Let $A \subseteq 2^\omega$ be a set of measure 1 which is
disjoint with $\bigcap_{\xi < \kappa} \F_\xi$. If $\F_\xi$'s are
ultrafilters then it follows that $A^\star \subseteq \bigcup_{\xi <
\kappa} \F_\xi$.
$\QED$

Let $\left\{I_n : n \in \omega\right\}$ be a partition of $\omega$ into finite
sets such that $|I_n| \geq 2^n$ for $n \in \omega$.
Define
$$A = \left\{X \subseteq \omega :  \exists a >1 \
\forall^\infty n \ \frac{|X \cap I_n|}
{|I_n|} > \frac{a}{3}\right\}  .$$
It is not hard to see that $\mu(A) = 1$.

Let $\F$ be an ultrafilter on $\omega$. Define a notion of forcing
${\cal Q}_\F$ as follows.
$${\cal Q}_\F = \left\{q \in [A]^{<\omega} : \exists a > 1 \
\left\{n : \frac{|I_n \cap \bigcap q|}{|I_n|} > \frac{a}{3^{|q|}} \right\} \in \F\right\}  .$$
Elements of ${\cal Q}_\F$ are ordered by inclusion.
\begin{lemma}\label{lem}
${\cal Q}_\F$ is powerfully ccc.
\end{lemma}
\Proof
We have to show that ${\cal Q}^n_\F$ satisfies countable chain
condition for every natural number $n$.

Let's start with the case when $n=1$.

 Suppose that ${\cal A} \subseteq {\cal Q}_\F$ is an uncountable
subset. Find $k \in \omega$, $a>1$ and  an uncountable set
${\cal A}' \subseteq {\cal A}$ such
that $|q|=k$  and 
$$ \left\{n : \frac{|I_n \cap \bigcap q|}{|I_n|} \geq \frac{a}{3^k} \right\} \in \F \ $$
for all $q \in {\cal A}'$.

We show that any sufficiently big subset of ${\cal A}'$
contains two compatible conditions.
We will use the following general observation.
\begin{lemma}
Let $(X, \nu)$ be a measure space with probability measure $\nu$.
Suppose that $1<b<a$ and $\varepsilon>0$. There exists a number $l$
such that if $A_1, \ldots, A_l$ are subsets of $X$ of measure $\geq a
\cdot \varepsilon$ then there are $i \neq j$ such that
$$\nu(A_i \cap A_j) \geq b^2 \cdot \varepsilon^2.$$
\end{lemma}
\Proof
Let $A_1, \ldots, A_l$ be sets of measure $a \cdot \varepsilon$.
Consider random variables $X_i$ given by characteristic function of
$A_i$ for $i \leq l$. Note that $X_i^2=X_i$ for $i \leq l$.

Suppose that $\nu(A_i \cap A_j) < b^2 \cdot \varepsilon^2$ for $i
\neq j$. In particular $E(X_i X_j) < b^2 \cdot \varepsilon^2$ for $i
\neq j$. Recall that ${\bf E}(X^2) \geq \lft1({\bf E}(X)\rgt1)^2$ for every
random variable $X$.

We compute
$$\eqalign{
0&\leq{\bf E}\lft2(\lft1((\sum_{i=1}^l X_i)-la\varepsilon\rgt1)^2\rgt2)\cr
&={\bf E}\lft1(\sum_{i=1}^l X_i^2+\sum_{i\neq j}X_iX_j
   -2la\varepsilon\sum_{i=1}^l X_i+(la\varepsilon)^2\rgt1)\cr
&\leq la\varepsilon+l(l-1)b^2\varepsilon^2
    -2(la\varepsilon)^2+(la\varepsilon)^2\cr
&=l\lft1(a\varepsilon-b^2\varepsilon^2+l(b^2-a^2)\varepsilon^2\rgt1)
.}$$
As $b<a$, the last line here is negative for large l, a contradiction.~$\QED$

Let $\nu_n$ be the uniform measure on $I_n$ i.e. $\nu(A) = |A| \cdot
|I_n|^{-1}$ for $A \subseteq I_n$. Fix $1<b<a$ and let $l$ be a number
from the above lemma chosen for $\varepsilon=3^{-k}$.

Let $q_1, \ldots, q_l \in {\cal A}'$.
Find $Y \in \F$ such that for all $i \leq l$ and $n \in Y$
$$\frac{|I_n \cap \bigcap q_i|}{|I_n|} \geq
\frac{a}{3^k} \ .$$
By the lemma for every $n \in Y$ there exist $i \neq j$ such that
$$\frac{|I_n \cap \bigcap (q_i \cup q_j)|}{|I_n|} \geq
\frac{b}{9^k} \ .$$
Since $\F$ is an ultrafilter there exist $i \neq j$  such that
$$\left\{n : \frac{|I_n \cap \bigcap (q_i \cup q_j)|}{|I_n|} \geq
\frac{b}{9^k}\right\} \in{\cal F} \ .$$
It follows that $q_i \cup q_j \in {\cal Q}_\F$. 

Note that in fact we proved that for every $m$ there exists $l$ such
that for every subset $X \subseteq {\cal A}'$ of size $l$ there are
$q_1, \ldots, q_m \in X$ such that $q_i \cup  q_j \in {\cal
Q}_\F$ for $i,j \leq m$.

Suppose that $n>1$. For simplicity assume that $n=2$, the general case
is similar.

Let ${\cal A}$ be an uncountable subset of ${\cal Q}_\F^2$.
Without loss of generality we can assume there are numbers $k_1$,
$k_2$ and $a$ such that
for every $\<q^1,q^2\> \in {\cal A}$, $|q^1|=k_1$, $|q^2|=k_2$ and
$$ \left\{n : \frac{|I_n \cap \bigcap q^i|}{|I_n|} \geq
\frac{a}{3^{k_i}} \right\} \in \F  \hbox{ for } i=0,1 . $$
To get two elements of ${\cal A}$ which are compatible first apply the
above remark to get a large subset $X \subseteq {\cal A}$ with first
coordinates being pairwise compatible and then apply the case $n=1$ to
the second coordinates of conditions in $X$.~$\QED$

Notice that if $G$ is a ${\cal Q}_\F$-generic filter over $\V$ then
the set
$$\left\{\bigcap q : q \in G \right\} \cup \F_0$$
generates a non-principal filter on $\omega$.

Let ${\cal P}_\F = \lim_{n \rightarrow \infty} {\cal Q}_\F^n$ be a finite
support product of countably many copies of ${\cal Q}_\F$.
By \ref{lem},  ${\cal P}_\F$ satisfies the countable chain
condition.

\begin{lemma}\label{crucial}
$\forces_{{\cal P}_\F} A \cap \V \hbox{ is the union of countably many
filters} .$
\end{lemma}
\Proof
Let $G$ be a ${\cal P}_\F$-generic filter over $\V$.

For $n \in \omega$,
let $G_n = \left\{ p(n): \dom(p)=\left\{n\right\}, p \in G\right\}$.
Since $G_n$ is a ${\cal
Q}_\F$-generic filter let $\F_n$ be a filter on $\omega$ generated by
$G_n$ as above.

We show that $\V[G] \models A \cap \V \subseteq \bigcup_{n \in \omega}
\F_n$.

Suppose that $\bar{p} \in {\cal P}_\F$ and $X \in A$. Find $n \in
\omega$ such that $n \not \in \dom(\bar{p})$. Let $\bar{q}=\bar{p}
\cup \<n,\left\{X\right\}\>$. It is clear that $\bar{q} \forces X \in \F_n$.
$\QED$

Now we can finish the proof of \ref{main}.

Let $\V \models 2^{\aleph_{0}}>\aleph_1$. Let $\<{\cal P}_\xi, \dot{{\cal
Q}}_\xi : \xi < \omega_1\>$ be a finite support iteration such that
$$\forces_\xi \dot{{\cal Q}}_\xi \simeq {\cal P}_{\F_\xi} \hbox{ for
some ultrafilter } \F_\xi\ .$$
Let ${\cal P}_{\omega_1} = \lim_{\xi < \omega_1} {\cal P}_\xi$.

Let $G$ be a ${\cal P}_{\omega_1}$-generic filter over $\V$. We will
show that $\V[G]$ is the model we are looking for.

Let $\left\{\F^\xi_n : \xi < \omega_1, n \in \omega\right\}$ be the family of
filters added by $G$.
Without loss of generality we can assume that they are ultrafilters.
To finish the proof it is enough to show that
$$A \subseteq \bigcup_{\xi < \omega_1} \bigcup_{n \in \omega} \F^\xi_n
\ .$$
Suppose that $X \in A$. There exists $\xi < \omega_1$ such that $X \in
\V[G \cap {\cal P}_\xi]$. By \ref{crucial} there exists $n \in
\omega$ such that $X \in F_n^{\xi +1}$. $\QED$

\end{document}